\documentclass[12pt]{article}
\usepackage{amsmath,amssymb}
\baselineskip 7pt \lineskip 7pt \numberwithin{equation}{section}

\def \N{\hbox{$I\hskip -3pt N$}}
\def \Z{\hbox{$Z\hskip -5.2pt Z$}}

\def \C{\hbox{$C\hskip -5pt \vrule height 6pt depth 0pt \hskip 6pt$}}

\def\qed{\ \ \ifhmode\unskip\nobreak\fi\ifmmode\ifinner
         \else\hskip5pt\fi\fi
 \hbox{\hskip5pt\vrule width4pt height6pt depth1.5pt\hskip 1 pt}}

\def\cen{\centerline}

\def\z{\mathbb{Z}}
\def\c{\mathbb{C}}
\def\dim{\hbox{dim}}

\newfont{\df}{eufm10}

\def\ot{\otimes}

\def\dim{\hbox{\rm dim}\,}

\def\ot{\otimes}

\def\l{\lambda}

\def\Vir{\mbox{Vir}}

\def\cl{\centerline}

\def\vs{\vspace*}

\def\C{\mathbb{C}}

\def\Z{\mathbb{Z}}
\def\N{\mathbb{N}}

\newtheorem{theo}{Theorem}[section]
\newtheorem{lemm}[theo]{Lemma}
\newtheorem{coro}[theo]{Corollary}

\newtheorem{defi}[theo]{Definition}

\title{\bf  Classification of irreducible weight modules over $W$-algebra $W(2,2)$
\thanks{E-mail: liudong@hutc.zj.cn}}
\author{Dong Liu, \quad Gao Shoulan\\ Department of Mathematics, Huzhou Teachers College\\ Zhejiang Huzhou, 313000, China\\
Linsheng Zhu\\ Department of Mathematics, Changshu Institute of
Technology\\ Jiangsu Changshu, 215500, China}
\date{ }
\begin{document}
\maketitle

{\small {\bf Abstract. } We show that the support of an
irreducible weight module over the $W$-algebra $W(2, 2)$, which
has an infinite dimensional weight space, coincides with the
weight lattice and that all nontrivial weight spaces of such a
module are infinite dimensional. As a corollary, we obtain that
every irreducible weight module over the the $W$-algebra $W(2,
2)$, having a nontrivial finite dimensional weight space, is a
Harish-Chandra module (and hence is either an irreducible highest
or lowest weight module or an irreducible module of the
intermediate series).
 \vs{2pt}\par {\bf Key
Words:} The the $W$-algebra $W(2, 2)$, weight modules, support
\vs{2pt}\par{\it  Mathematics Subject Classification (2000)}:
17B56; 17B68.}

\vs{10pt}\par \cl{\bf1. \ Introduction}
\setcounter{section}{1}\setcounter{theo}{0}

The $W$-algebra $W(2, 2)$ was introduced in \cite{DZ} for the
study the classification of vertex operator algebras generated by
weight 2 vectors.

\begin{defi} The {\bf $W$-algebra} ${{\cal L}}=W(2, 2)$ is
a Lie algebra over $\c$ (the  field of complex numbers) with the
basis
$$\{x_n,I(n),C, C_1 | n \in \z\}$$
and the Lie bracket given by
$$[x_n,x_m]=(m-n)x_{n+m}+\delta_{n,-m}\frac{n^3-n}{12}C,\eqno (1.1)$$
$$[x_n,I(m)]=(m-n)I(n+m)+\delta_{n,-m}\frac{n^3-n}{12}C_1, \eqno (1.2)$$
$$[I(n),I(m)]=0, \eqno (1.3)$$
$$[{\cal L},C]=[{\cal L},C_1]=0. \eqno (1.4)$$
\end{defi}

The $W$-algebra $W(2,2)$ can be realized from the semi-product of
the Virasoro algebra Vir and the Vir-module ${\cal A}_{0, -1}$ of
the intermediate series in \cite{OR}. In fact, let $W=\c\{x_m\mid
m\in\z\}$ be the Witt algebra (non-central Virasoro algebra) and
$V=\c\{I(m)\mid n\in\z\}$ be a $W$-module with the action
$x_m\cdot I(n)=(n-m)I(m+n)$, then $W(2, 2)$ is just the universal
central extension of the Lie algebra $W\ltimes V$ (see \cite{OR}
and \cite{GJP}).  The $W$-algebra $W(2,2)$ studied in \cite{DZ} is
the restriction for $C_1=C$ of $W(2,2)$ in our paper.

The $W$-algebra $W(2,2)$ can be also realized from the so-called
{\it loop-Virasoro algebra} (see \cite{GLZ}). Let $\c[t, t^{-1}]$
be the Laurents polynomial ring over $\c$, then the loop-Virasoro
algebra $\tilde{VL}$ is the universal central extension of the
loop algebra $\Vir\ot\c[t^1, t^{-1}]$ and $W(2,
2)=\tilde{VL}/\c[t^2]$.

 The $W$-algebra
$W(2,2)$ is an extension of the Virasoro algebra and is similar to
the twisted Heisenberg-Virasoro algebra (see \cite{ADKP}).
However, unlike the case of the later, the action of $I(0)$ in
$W(2, 2)$ is not simisimple, so its representation theory is very
different from that of the twisted Heisenberg-Virasoro algebra in
a fundamental way.

Next we recall the definitions of $\z$-graded
${\cal L}$-modules.  If ${\cal L}$-module $V=\oplus
_{ m\in{\z}}V_{m}$ satisfies
$$
{\cal L}_m\cdot V_{n}\subset V_{ m+n},\quad\forall m,n\in{\z},
\eqno{(3.2)}
$$
then $V$ is called a {\it ${\z}$-graded ${\cal L}$-module}, and
$V_{m}$ is called a {\it homogeneous subspace of $V$ with degree}
$m\in{\z}$.

A ${\z}$-graded module $V$ is called {\it quasi-finite} if all
homogeneous subspaces are finite dimensional; a {\it uniformly
bounded module} if there exists a number $n\in{\bf N}$ such that
all dimensions of the homogeneous subspaces are $\le n$; a {\it
module of the intermediate series} if $n=1$.

For any ${\cal L}$-module $V$ and $\lambda\in \c$, set
$V_\lambda:=\bigl\{v\in V\bigm|x_0v=\lambda v\bigr\}$, which we
generally call the weight space of $V$ corresponding the weight
$\lambda$.

\par
An ${\cal L}$-module $V$ is called a weight module if $V$ is the
sum of all its weight spaces. For a weight module $V$ we define
$$\hbox{Supp}(V):=\bigl\{\lambda\in \c \bigm|V_\lambda\neq
0\bigr\},$$ which is generally called the weight set (or the
support) of $V$.

\par
 A nontrivial weight ${\cal L}$-module V is called a {weight module of intermediate
 series} if V is indecomposable and any weight spaces of V is one dimensional.
\par
A  weight ${\cal L}$-module V is called a {highest} (resp.
{lowest) weight module} with {highest weight} (resp. {highest
weight}) $\lambda\in \c$, if there exists a nonzero weight vector
$v \in V_\lambda$ such that

 \vskip 5pt 1) $V$ is generated by $v$ as  ${\cal L}$-module;

2) ${\cal L}_+ v=0 $ (resp. ${\cal L}_- v=0 $).
 \vskip 5pt

\noindent{\bf Remark.} For a highest (lowest) vector $v$ we always
suppose that $I_0v=c_0v$ for some $c_0\in \c$ although the action
of $I_0$ is not semisimple.

Obviously, if $M$ is an irreducible weight $\L$-module, then there
exists $\lambda\in\C$ such that ${\rm Supp}(M)\subset\lambda+\Z$.
So $M$ is a $\z$-graded module.

If, in addition, all weight spaces $M_\l$ of a weight $\L$-module
$M$ are finite dimensional, the module is called a {\it
Harish-Chandra module}. Clearly a highest (lowest) weight module
is a Harish-Chandra module.

Let $U:=U(\cal L)$ be the universal enveloping algebra of  ${\cal
L}$. For any $\lambda,c$ $\in \c$, let $I(\lambda,c, c_0, c_1)$ be
the left ideal of $U$ generated by the elements $$
\bigl\{x_i,I(i)\bigm|i\in \N \bigr\}\bigcup\bigl\{x_0-\lambda
\cdot 1, C-c\cdot 1, I_0-c_0\cdot 1, C_1-c_1\cdot1\bigr\}. $$ Then
the Verma module with the highest weight $\lambda$ over ${\cal L}$
is defined as
$$M(\lambda, c, c_0, c_1):=U/I(\lambda, c, c_0, c_1).$$
It is clear that $M(\lambda, c, c_0, c_1)$ is a highest weight
module over ${\cal L}$ and   contains a unique maximal submodule.
Let $V(\lambda, c, c_0, c_1)$ be the unique irreducible quotient
of $M(\lambda, c, c_0, c_1)$.

The following result was given in \cite{DZ}.
\begin{theo} \cite{DZ} The Verma module $M(\lambda, c, c_0, c_1)$ is irreducible if and only if
$\frac{m^2-1}{12}c_1+2c_0\ne 0$ for any nonzero integer $m$.
\end{theo}

The classification of Harish-Chandra modules over the W-algebra
$W(2,2)$ was given in \cite{LLZ}.
\begin{theo}\label{T2}\cite{LLZ} A Harish-Chandra module
$\cal L$-module $V$ is a highest weight module or lowest weight
module or a module of the intermediate series.
\end{theo}

An irreducible weight module $M$ is called a {\it pointed module}
if there exists a weight $\l\in\C$ such that ${\rm dim\,}V_\l=1$.
Xu posted the following in \cite{X}:

\noindent{\bf Problem 1.1~}~{\it Is any irreducible pointed module
over the Virasoro algebra a Harish-chandra module?}\vs{4pt}\par An
irreducible weight module $M$ is called a {\it mixed  module} if
there exist $\l\in\C$ and $i\in\Z$ such that ${\rm
dim\,}V_\l=\infty$ and ${\rm dim\,}V_{\l+i}<\infty$. The following
conjecture was posted in \cite{M}:

\noindent{\bf Conjecture 1.2~}~{\it There are no irreducible mixed
module over the Virasoro algebra.}\vs{4pt}\par Mazorchuk and Zhao
\cite{MZ} gave the positive answers to the above question and
conjecture to the Virasoro algebra, Shen and Su \cite{SS} also
also gave a similar result for the twisted Heisenberg-Virasoro
algebra.

 In this paper, we also give the positive answers to the above
 question and conjecture for the $W$-algebra $W(2, 2)$.
 Due to many differences between the $W(2,2)$ and the twisted Heisenberg-Virasoro
algebra, some new methods are given in our paper.

Our main result is the following:
\begin{theo}\label{tmain} Let M be an irreducible
weight $\L$-module. Assume that there exists $\lambda\in\C$ such
that ${\rm dim\,}M_\lambda=\infty$. Then ${\rm
Supp}(M)=\lambda+\Z$, and for every $k\in\Z$, we have ${\rm
dim\,}M_{\lambda+k}=\infty$.
\end{theo}

The paper is organized as follow: Some lemma for the proof of
Theorem~\ref{tmain} are given in Section 2. The Proof of the main
Theorem is given in Section 3 where some corollaries from this
theorem are also discussed.

\section{Point modules over the $W$-algebra}
\setcounter{equation}{0} \vs{1pt}\par We first recall a main
result about the weight Virasoro-module in [MZ]:

\begin{theo} Let $V$ be an irreducible weight
Virasoro-module. Assume that there exists $\lambda\in \C$, such
that ${\rm dim\,}V_\lambda=\infty$. Then ${\rm
Supp}(V)=\lambda+\Z$, and for every $k\in\Z$, we have ${\rm
dim\,}V_{\lambda+k}=\infty$.
\end{theo}

\begin{lemm}\label{HH}  Assume that there exists $\mu\in\C$ and
a non-zero element $v\in M_\mu$, such that
$$I_1v=L_1v=L_{-1}I_2v=L_2v=0\mbox{ \ \  or \ \ }I_{-1}v=L_{-1}v=L_1I_{-2}v=L_{-2}v=0.$$ Then $M$ is a
Harish-Chandra module.
\end{lemm}

\vs{4pt}\par{\it Proof.~} Suppose that $I_1v=L_1v=L_2v=0$ for
$v\in V_\mu$, it is clear that $L_{>0}v=0$ and $I_mv=0$ for $m\ge
3$. Moreover $L_{>0}I_2v=0$ and $I_mI_2v=0$ for $m\ge 3$ or $m=1$.

But $L_{-1}I_2v=0$, then $L_1L_{-1}I_2v=[L_{1},
L_{-1}]I_2v+L_{-1}L_{1}I_2v=-{1\over2}L_0I_2v=0$. So $I_2v=0$ if
$\mu\ne -2$. Then ${\cal L}_{>0}v=0$. Hence $v$ is a highest
weight vector, and hence, $M$ is a Harish-Chandra module.

If $\mu=-2$ and $w=I_2v\ne 0$, then $L_nw=L_nI_2v=[L_n,
I_2]v+I_2L_nv=0$ for any $n\in\N$. Moreover $I_1w=0$ and
$L_{-1}I_2w=[L_{-1}, I_2]I_2v+I_2L_{-1}I_2v=0$. Then $I_2w=0$
since $I_2w\notin V_{0}$. So ${\cal L}_{>0}w=0$. Hence $w$ is
either a highest weight vector, and hence, $M$ is a Harish-Chandra
module. Similar for the lowest weight case.$\hfill\Box$

Assume now that $M$ is an irreducible weight $\L$-module such that
there exists $\lambda\in\C$ satisfying ${\rm
dim\,}M_\lambda=\infty$.

\begin{lemm}\label{l1} There exists at most one $i\in\Z$ such
that ${\rm dim\,}M_{\lambda+i}<\infty$.
\end{lemm}
\noindent{\bf Proof.~}~Assume that $${\rm
dim\,}M_{\lambda+i}<\infty\mbox{ \ \  and  \ }{\rm
dim\,}M_{\lambda+j}<\infty\mbox{ \ \ \ for some different \ \
}i,j\in\Z.$$ Without loss of generality, we may assume $i=1$ and
$j>1$. Set
$$
\begin{array}{ll}
V:=&\mbox{Ker}(I_1:M_\lambda\rightarrow
M_{\lambda+1})\cap\mbox{Ker} (L_1, L_{-1}I_2:M_\lambda\rightarrow
M_{\lambda+1})\cap\mbox{Ker}(I_j:M_\lambda\rightarrow
M_{\lambda+j})\\[7pt]&
\cap\,\mbox{Ker}(L_j:M_\lambda\rightarrow
M_{\lambda+j}),\end{array}$$ which is a subspace of $M_\l$. Since
$${\rm dim\,}M_\lambda=\infty,\ \ \ {\rm dim\,}M_{\lambda+1}<\infty\mbox{ \ \
and \ \ }{\rm dim\,}M_{\lambda+j}<\infty,$$ we have, ${\rm
dim\,}V=\infty.$ Since $$[L_1,L_k]=(k-1)L_{k+1}\neq 0\mbox{ \ \
and \ }[I_1,L_l]=(l-1)I_{l+1}\neq 0\mbox{ \ \ for \ \ }k, l\in\Z,\
\ k, l\ge 2,$$ we get
$$
\begin{array}{ll}
L_kV=0,& k=1,j,j+1,j+2,\cdots,\ \mbox{\ \ and \ \ }\\[7pt]
I_lV=0,& l=1,j,j+1, j+2\cdots.\end{array}\eqno(2.1)$$If there
would exist $0\neq v\in V$ such that $L_2v=0$, then
$I_1v=L_1v=L_{-1}I_2v=L_2v=0$ and $M$ would be a Harish-Chandra
module by Lemma 2.2. It is a contradiction. Hence $L_2v\neq 0$ for
all $v\in V$. In particular,
$${\rm dim\,}L_2V=\infty.$$ Since ${\rm dim\,}M_{\lambda+1}<\infty$,
and the actions of $I_{-1}$ and $L_{-1}$ on $L_2V$ map $L_2V$
(which is an infinite dimensional subspace of $M_{\l+2}$) to
$M_{\l+1}$ (which is finite dimensional), there exists $0\neq w\in
L_2V$ such that $I_{-1}w=L_{-1}w=0$. Let $w=L_2v$ for some $v\in
V$. For all $k\geq j$, using (2.1), we have
$$L_kw=L_kL_2v=L_2L_kv+(2-k)L_{k+2}v=0+0=0.$$ Hence $L_kw=0$ for
all $k=1,j,j+1,j+2,\cdots.$ Since
$$[L_{-1},L_l]=(l+1)L_{l-1}\neq 0\ \mbox{ \ and \ } [I_{-1},I_l]=(l+1)I_{l-1}\neq 0\ \mbox{\ for
all \ }l>1,$$ we get inductively $L_kw=I_kw=0$ for all
$k=1,2,\cdots.$ Hence $M$ is a Harish-Chandra module by Lemma 2.2. A
contradiction. The lemma follows.$\hfill\Box$

Because of Lemma 2.3, we can now fix the following notation: $M$ is
an irreducible weight $\L$-module, $\mu\in\C$ is such that ${\rm
dim\,}M_\mu<\infty$ and ${\rm dim\,}M_{\mu+i}=\infty$ for every
$i\in\Z\setminus\{0\}$.

\begin{lemm}\label{ll2}

Let $0\neq v\in M_{\mu-1}$ and $\mu\ne -1$ such that
$I_1v=L_1v=L_{-1}I_2v=0$. Then

(1) There exists a nonzero $u\in M$ such that $L_1v=I_mv=0$ for
all $m\ge1$.

(2) $I_mL_2v=0$ for all $m\ge1$.
\end{lemm}
\noindent{\bf Proof.~} Since $L_{-1}I_2v=0$, then
$L_1L_{-1}I_2v=[L_{1},
L_{-1}]I_2v+L_{-1}L_{1}I_2v=-{1\over2}L_0I_2v=0$. So $I_2v=0$
since $\mu\ne -1$. By $[L_1, I_k]=(k-1)I_{k+1}$ we have $I_kv=0$
for all $k\ge 2$. Moreover $I_mL_2v=[I_m,
L_2]v+L_2I_mv=(2-m)I_{m+2}v+L_2I_mv=0$. $\hfill\Box$

\begin{lemm}\label{ll3}

Let $0\neq w\in M_{\mu+1}$ and $\mu\ne 1$ such that
$I_{-1}w=L_{-1}w=L_{1}I_{-2}w=0$. Then

(1)$L_{-1}w=I_{-m}w=0$ for all $m\ge1$.

(2) $I_{-m}L_{-2}w=0$ for all $m\ge1$.

\end{lemm}

\noindent{\bf Proof.~} It is similar to that in Lemma 2.4.
$\hfill\Box$

\section{\bf Proof of Theorem \ref{tmain}}

\noindent{\bf Proof of Theorem \ref{tmain}.} Due to Lemma
\ref{l1}, we can suppose that $\dim M_\mu<+\infty$ and $\dim
M_{\mu+i}=+\infty$ for all $i\in\z, i\ne 0$.

Set
\begin{multline*}
V:=\mbox{Ker}\{L_1:M_{\mu-1}\rightarrow M_\mu\}\cap
\mbox{Ker}\{I_1:M_{\mu-1}\rightarrow M_\mu\}\\  \cap
\mbox{Ker}\{L_{-1}I_2 :M_{\mu-1}\rightarrow
M_\mu\}\cap\mbox{Ker}\{L_{-1}L_2 :M_{\mu-1}\rightarrow
M_\mu\}\subset M_{\mu-1}. \end{multline*}

For any $v\in V$, $L_1v=I_1v=L_{-1}I_2v=0$

Since ${\rm dim\,}M_{\mu-1}=\infty$ and ${\rm dim\,}M_\mu<\infty$,
we have ${\rm dim\,}V=\infty$. For any $v\in V$, consider the
element $L_2v$. By Lemma 2.2, $L_2v=0$ would imply that $M$ is a
Harish-Chandra module, a contradiction. Hence $L_2v\neq 0$, in
particular, ${\rm dim\,}L_2V=\infty$.

Since the actions of $I_{-1}$, $L_{-1}$,$L_{1}L_{-2}$ and
$L_{1}I_{-2}$ on $L_2V$ map $L_2V$ (which is an infinite
dimensional subspace of $M_{\mu+1}$) to $M_\mu$ (which is finite
dimensional), there exists $w=L_2v\in L_2V$ for some $v\in V$,
such that $0\ne w\in M_{\mu+1}$ and
$I_{-1}w=L_{-1}w=L_{1}I_{-2}w=L_{1}L_{-2}w=0$.

(1) If $\mu\ne \pm1$, then $$I_kw=0,\ k=1,2,\cdots\eqno(3.1)$$
from Lemma \ref{ll2} and
$$I_{-k}w=0,\ k=1,2,\cdots\eqno(3.2)$$ from Lemma \ref{ll3}.

 This means that $I_k$ act trivially on $M$ for all $k\in\Z$, and
so $M$ is simply an irreducible module over the Virasoro algebra.
Thus, Theorem 1.3 follows from Theorem 2.1 in the case
$\mu\ne\pm1$.

(2) If $\mu=\pm1$, we only show that $\mu=1$ is not possible and
for $\mu=-1$ the statement will follow by applying the canonical
involution on ${\cal L}$.

In fact, if $\mu=1$, then for $v\in V$,
$L_1v=I_1v=L_{-1}I_2v=L_{-1}L_2v=L_0v=0$. By Lemma \ref{ll2}, we
have $I_kv=0, k=1, 2, \cdots$.

For any $v\in V$, consider the element $L_2v$. By Lemma 2.2,
$L_2v=0$ would imply that $M$ is a Harish-Chandra module, a
contradiction. Hence $L_2v\neq 0$, in particular, ${\rm
dim\,}L_2V=\infty$.

Since the actions of $I_{-1}$, $L_{-1}$ and $L_{1}I_{-2}$ on
$L_2V$ map $L_2V$ (which is an infinite dimensional subspace of
$M_{2}$) to $M_1$ (which is finite dimensional), there exists
$w=L_2v\in L_2V$ for some $v\in V$, such that $w\neq 0$ and
$I_{-1}w=L_{-1}w=L_{1}I_{-2}w=0$. Moreover we have $$I_kw=0, k=1,
2,\cdots\eqno(3.3)$$ from Lemma \ref{ll2}. So $$I_0w=0.
\eqno(3.4)$$

If $L_1w=L_1L_2v=0$, then from $L_{-1}L_2v=0$ we have
$L_{1}L_{-1}L_2v=[L_1, L_{-1}]L_2v+L_{-1}L_1L_2v=0$. So
$L_0L_2v=2L_2V=0$ since $L_2v\in M_2$. Hence $L_2v=0$ and then $M$
is a highest weight module.

Then we can suppose that $L_1w\ne 0$ for any $w\in L_2V$.

 For any $w\in L_2V$, consider the element $L_{-2}w$. If
$L_{-2}w=0$, then $L_{-k}w=I_{-k}w=0,\ k=1,2,\cdots$. Then M is a
Harish-Chandra module.  Hence $L_{-2}L_2V\ne 0$, in particular,
${\rm dim\,}L_{-2}L_2V=\infty$. Let $W=L_2V$, then $L_1$ maps
$L_{-2}W$ to $M_1$ has infinite dimensional Kernel $K$. Let $0\ne
L_{-2}w\in K$, then $L_1L_{-2}w=0$. But $L_1L_{-2}=[L_1,
L_{-2}]+L_{-2}L_1$ and $[L_1, L_{-2}]w=(-3)L_{-1}w=0$, hence
$L_{-2}L_1w=0$. Setting $u=L_1w\ne 0$, we have $L_{-2}u=0$,
$I_{-1}u=I_{-1}L_1w=[I_{-1}, L_1]w+L_{1}I_{-1}w=0$. Moreover by
induction we have $I_{-m}u=0$ for all $m\ge 3$.

 $I_mu=[I_m,
L_1]w+L_1I_mw=(1-m)I_{m+1}w+L_1I_mw=0$ for all $m\ge 0$ by (3.3)
and (3.4). So $I_2u={1\over6}[L_{-2}, I_4]u=0$. Then $I_ku=0$ for
all $k\in\z$.

By $L_{-2}u=0$, we have $I_2L_{-2}u=0$. Therefore $c_1=0$.

This means that $I_k, k\in\z, C_1$ act trivially on the
irreducible $M$ for all $k\in\Z$, and so $M$ is simply a module
over the Virasoro algebra. Thus, Theorem 1.3 follows from Theorem
2.1. $\hfill\Box$

Theorem \ref{tmain} also implies the following classification of
all irreducible weight $\L$-modules which admit a nontrivial
finite dimensional weight space:
\begin{coro} Let $M$ be an irreducible weight
$\L$-module. Assume that there exists $\lambda\in\C$ such that
$0<{\rm dim\,}M_\lambda<\infty$. Then $M$ is a Harish-Chandra
module. Consequently, $M$ is either an irreducible highest or
lowest weight module or an irreducible module from the
intermidiate series. \end{coro}

\noindent{\bf Proof.}~Assume that $M$ is not a Harish-Chandra
module. Then there should exists $i\in\Z$ such that ${\rm
dim\,}M_{\lambda+i}=\infty$. In this case, Theorem 1.3 implies
${\rm dim\,}M_\lambda=\infty$, a contradiction. Hence $M$ is a
Harish-Chandra module, and the rest of the statement follows from
Theorem \ref{T2}.$\hfill\Box$
\par

\vskip30pt \centerline{\bf ACKNOWLEDGMENTS}

\vskip15pt Project is supported by the NNSF (Grant 10671027,
10701019, 10571119), the ZJZSF(Grant Y607136), and Qianjiang
Excellence Project of Zhejiang Province (No. 2007R10031). Authors
give their thanks to Prof. Congying Dong for his useful comments.
\vskip30pt

\def\refname{\cen{\bf REFERENCES}}

\end{document}